# A Robust Mixed Integer Optimization Model to Utilize Regenerative Energy of Trains in a Railway Network

Shuvomoy Das Gupta, J. Kevin Tobin and Lacra Pavel[*][†]

**Abstract**

In this paper we present a robust mixed integer optimization model to utilize regenerative braking energy produced by trains in a railway network. An electric train produces regenerative energy during braking, which is often lost in present technology. To utilize this energy we calculate a timetable which maximizes the total overlapping time between the braking and accelerating phases of suitable train pairs, so that the regenerative energy of a braking train can be transferred to a suitable accelerating one. We apply our optimization model to different instances of a railway network for a time horizon spanning six hours. For each instance, our model finds an optimal timetable very quickly (largest running time being 86.64 seconds). Compared to the existing timetables, we observe significant increase in utilization of regenerative energy for every instance.

***Index Terms***— Railway networks, energy efficiency, regenerative braking, train scheduling, mathematical optimization.

# 1 Introduction

## 1.1 Background and motivation

Application of mathematical optimization techniques to efficient energy management of electric vehicles has gained lot of interest in the past few years[1, 2, 3, 4]. Among all of the public transport modes, trains are often preferred by the passengers for providing higher capacity and safety. In most railway networks, trains use electricity as their primary source of energy and are equipped with regenerative braking technology. When a train makes a trip from an origin platform to a destination platform, its optimal speed profile consists of four phases: 1) accelerating, 2) speed holding, 3) coasting and 4) braking [5], as shown in Figure 1. in a qualitative manner. Most of the energy required by the train is consumed during the accelerating phase. During the speed holding phase the energy consumption is negligible compared to accelerating phase, and during the coasting phase there is no need for energy. When the train brakes, it produces energy by using a regenerative braking mechanism. This energy is called regenerative braking energy. Naturally, proper utilization of regenerative braking energy of trains can lead to significant energy savings.

[*]This work was supported by NSERC-CRD and Thales Canada Inc.

[†]S. D. Gupta and L. Pavel are with the Edward S. Rogers Department of Electrical and Computer Engineering, University of Toronto, and J. K. Tobin is with Thales Canada Inc. e-mail: (shuvomoy.dasgupta@mail.utoronto.ca, Kevin.Tobin@thalesgroup.com, pavel@control.utoronto.ca)

However, transfer of the regenerative braking energy back to the electrical grid requires specialized technology such as reversible electrical substations [6, page 30]. On the other hand, storing it using current technology is not viable, as the storage options e.g., via supercapacitors, flywheels [7] are very expensive and have high discharge rates [8, page 66], [9, page 92]. A better strategy to utilize the regenerative braking energy of a train, that can be implemented with current technology, is to synchronize its braking phase with the accelerating phase of another nearby train operating under the same electrical substation. A positive overlapping time that arises from such a synchronization process enables transfer of the regenerative braking energy of the first train to the second one via the overhead contact line or a third rail [10], and can save the electrical energy that would be lost otherwise. Our objective is to design a railway timetable that contains the arrival and departure time of every train to and from all the platforms it visits, such that the duration of the **s**ynchronization **p**rocesses between **s**uitable **t**rain **p**airs (**SPSTP**s) is maximized subject to the different constraints present.

## 1.2 Related work

Over the past three decades, the general timetabling problem in a railway network has been studied extensively with objectives such as total cost, total passenger delay, quality of the service *etc.*[11]. However, a few among them attempt to determine a timetable to utilize regenerative braking energy of trains; all of them have appeared in the last three years.

The work by Peña-Alcaraz et al. [12] formulates a Mixed Integer Programming (**MIP**) model. In this model the objective is to maximize the total duration of all possible synchronization processes between all train pairs. The model is then successfully applied to line 3 of the Madrid underground system in Spain. However, for many railway networks, inclusion of all the train pairs in the objective may not be realistic, because energy transfer between trains that are not close to each other suffers from drastic transmission loss. Besides, inclusion of all such train pairs can make the optimization problem computationally intractable even for a medium sized networks, as each such train pair in the model corresponds to one binary variable. Because of the absence of connection and turn-around constraints, the model is only applicable to a single train-line.

A more tractable MIP model is presented in [13]. The model considers only the suitable train pairs for transfer of regenerative energy and can be applied to railway network of arbitrary topology. The optimization model is then applied to the Dockland Light Railway. The results show a significant increase in the total duration of the synchronization process. Although such increase, in principle, may increase the total savings in regenerative energy, the actual savings in terms of energy is not calculated. Another drawback of this model is that it lacks robustness, as it assumes that even if the trip time changes, the duration of the associated braking and accelerating stay the same, which is not the case in real life.

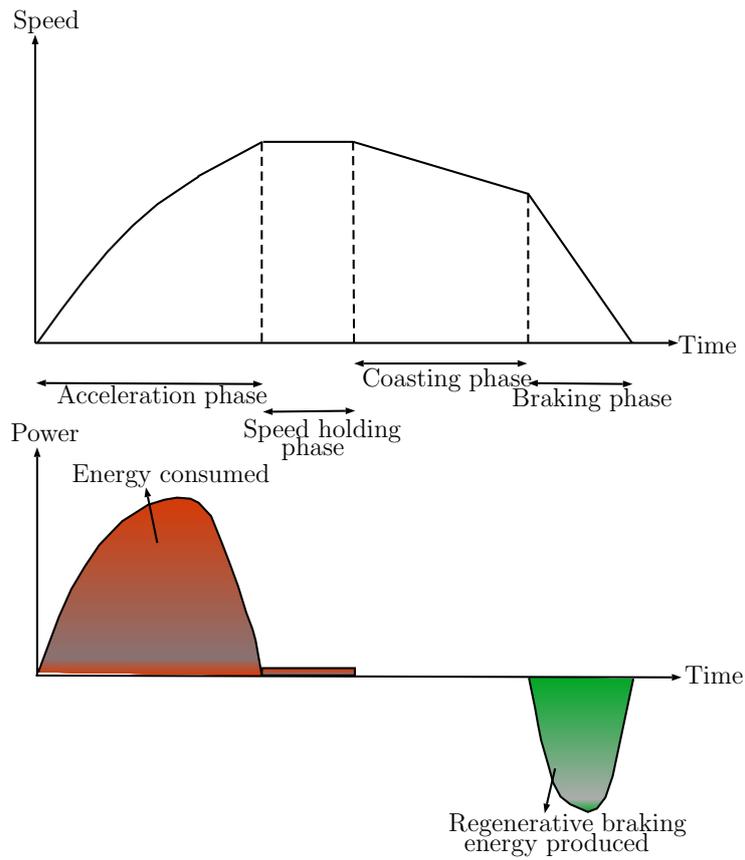

**Figure 1:** The speed profile of a train (top), and corresponding energy consumed and produced in accelerating phase and braking phase respectively (bottom).

Other relevant works implement meta-heuristics such as simulated annealing [14] and genetic algorithm [14] and [15]. These models cannot provide any guarantee of optimality for the generated timetable, because a meta-heuristic approach depends on manual tuning of the parameters.

An analytical study of a simplified and periodic railway schedule appears in [16] with focus on storage of regenerative energy. The model does not incorporate periodic event scheduling constraints [17] used to formulate periodicity in a timetable, so it is not possible to apply the model to some cases, *e.g.*, where trains associated with the second period enter the network while the trains associated with the first period are still running.

## 1.3 Contributions

Our contributions in this paper can be summarized as follows:

- Using mathematical programming techniques, we propose a robust MIP model to utilize the regenerative braking energy of trains in a railway network. In comparison with the existing works, our model can be applied to a railway network of arbitrary topology.

- For most of the existing railway networks, the railway management has a feasible timetable, which we exploit to devise an optimization strategy. This strategy produces a smaller search space and is computationally more tractable in comparison with the related works. We prove that all possible cases arising from an SPSTP are modelled accurately by our model.

- We apply our optimization model to different instances of a real-life railway network and find that our model gives optimal timetables very quickly. We observe significant increase in the utilization of regenerative braking energy.

## 1.4 Organization

This paper is organized as follows. In Section 2 we describe the notation and notions to be used in this paper. In Section 3 we define an SPSTP mathematically, propose our optimization strategy and formulate a mixed integer linear objective function. Then in Section 4 we model and justify various constraints present in the railway network. We propose the full optimization model in Section 5. In Section 6 we apply our model to realistic instances of different size. Section 7 presents the conclusion.

## 2 Notation and notions

All the sets described in this paper are strictly ordered and finite unless otherwise specified. The cardinality and the $i$th element of such a set $S$ is denoted by $|S|$ and $S(i)$ respectively.

Consider a railway network where the set of all stations is denoted by $\mathcal{S}$. The set of all platforms in the railway network is indicated by $\mathcal{N}$. A directed arc between two distinct and non-opposite platforms is called a track. The set of all tracks is represented by $\mathcal{A}$. The directed graph of the railway network is expressed by $\mathcal{G} = (\mathcal{N}, \mathcal{A})$. A train-line or line is a directed path, where the set of nodes represents non-opposite platforms and the set of arcs represents non-opposite tracks. A crossing-over is a special type of directed arc that connects two train-lines. If a train arrives at the terminal platform of a train-line, turns around by traversing the crossing-over, and starts travelling through another train-line, then the same physical train is treated and labelled functionally as two different trains by the railway management [18, page 41]. The set of all trains to be considered in our problem is denoted by $\mathcal{T}$. The sets of all platforms and all tracks visited by a train $t$ in chronological order are denoted by $\mathcal{N}^t \subseteq \mathcal{N}$ and $\mathcal{A}^t \subseteq \mathcal{A}$ respectively. The train-path of a train is the directed path that contains all platforms and tracks visited by it in chronological order.

The decision variables to be determined are the arrival and departure times of trains to and from the associated platforms respectively. Let $a_i^t$ be the arrival time of train $t \in \mathcal{T}$ at platform $i \in \mathcal{N}^t$ and $d_j^t$ be the departure time of train $t$ from platform $j \in \mathcal{N}^t$. Our objective is to maximize the duration of total overlapping time of all SPSTPs subject to the constraints, which would make it possible to save significant amount of electrical energy produced by the braking trains by transferring it to the accelerating trains.

## 3 Modelling the objective

At first, we need to characterize the train pairs and the associated platform pairs necessary to describe the SPSTPs. The platform pairs to consider are those opposite to each other and powered by the same electrical substations, because the transmission loss in transferring electrical energy between them is negligible. The set that contains all such platform pairs is denoted by $\Omega$. Consider any such platform pair $(i, j) \in \Omega$, and let $\mathcal{T}_i \subseteq \mathcal{T}$ be the set of all the trains which arrive at, dwell and then depart from platform $i$. Suppose, $t \in \mathcal{T}_i$. Now, we are interested to find another train $\tilde{t}$ on platform $j$, i.e., $\tilde{t} \in \mathcal{T}_j$, which along with $t$ would form a suitable pair for the transfer of regenerative braking energy. To achieve this, we start with an initial feasible timetable for the railway, which represents the desired service to be delivered. Too much deviation from it is not desired. For every train $t$, this timetable provides a feasible arrival time $\bar{a}_i^t$ and a feasible departure time $\bar{d}_i^t$ to and from every platform $i \in \mathcal{N}^t$ respectively. Intuitively, among all the trains that go through platform $j$, the one which is temporally closest to $t$ in the initial timetable would be the best candidate to form a pair with $t$. The temporal proximity can be of two types with respect to $t$, which results in the following definitions.

**Definition 1.** *Consider any $(i, j) \in \Omega$. For every train $t \in \mathcal{T}_i$, the train $\overrightarrow{\tilde{t}} \in \mathcal{T}_j$ is called the* ***temporally closest train to the right of*** $t$ *if*

$$\vec{t} = \underset{t' \in \{x \in \mathcal{T}_j : 0 \leq \frac{\bar{a}_j^x + \bar{d}_j^x}{2} - \frac{\bar{a}_i^t + \bar{d}_i^t}{2} \leq r\}}{\operatorname{argmin}} \left\{ \left| \frac{\bar{a}_i^t + \bar{d}_i^t}{2} - \frac{\bar{a}_j^{t'} + \bar{d}_j^{t'}}{2} \right| \right\}, \tag{1}$$

*where $r$ is an empirical parameter determined by the timetable designer and is much smaller than the time horizon of the entire timetable.*

**Definition 2.** *Consider any $(i,j) \in \Omega$. For every train $t \in \mathcal{T}_i$, the train $\overleftarrow{t} \in \mathcal{T}_j$ is called the **temporally closest train to the left of** $t$ if*

$$\overleftarrow{t} = \underset{t' \in \{x \in \mathcal{T}_j : 0 < \frac{\bar{a}_i^t + \bar{d}_i^t}{2} - \frac{\bar{a}_j^x + \bar{d}_j^x}{2} \leq r\}}{\operatorname{argmin}} \left\{ \left| \frac{\bar{a}_i^t + \bar{d}_i^t}{2} - \frac{\bar{a}_j^{t'} + \bar{d}_j^{t'}}{2} \right| \right\}. \tag{2}$$

**Definition 3.** *Consider any $(i,j) \in \Omega$. For every train $t \in \mathcal{T}_i$, the train $\tilde{t} \in \mathcal{T}_j$ is called the **temporally closest train to** $t$ if*

$$\tilde{t} = \underset{t' \in \{\vec{t}, \overleftarrow{t}\}}{\operatorname{argmin}} \left\{ \left| \frac{\bar{a}_i^t + \bar{d}_i^t}{2} - \frac{\bar{a}_j^{t'} + \bar{d}_j^{t'}}{2} \right| \right\}. \tag{3}$$

*If both $\vec{t}$ and $\overleftarrow{t}$ are temporally equidistant from $t$, we pick one of them arbitrarily.*

Any SPSTP can be described by specifying the corresponding $i$, $j$, $t$ and $\tilde{t}$ by using the definitions above. We construct a set of all the SPSTPs, which we denote by $\mathcal{E}$. Each element of this set is a tuple of the form $(i, j, t, \tilde{t})$. Because $\tilde{t}$ is unique for any $t$ in each element of $\mathcal{E}$, we can partition $\mathcal{E}$ into two sets denoted by $\vec{\mathcal{E}}$ and $\overleftarrow{\mathcal{E}}$, which contain elements of the form $(i, j, t, \vec{t})$ and $(i, j, t, \overleftarrow{t})$ respectively. For every $(i, j, t, \vec{t}) \in \vec{\mathcal{E}}$ our strategy is to synchronize the accelerating phase of $t$ with the braking phase of $\vec{t}$. On the other hand, for every $\overleftarrow{t} \in \overleftarrow{\mathcal{E}}$, then it would be convenient to synchronize the accelerating phase of $\overleftarrow{t}$ with the braking phase of $t$. For every $(i, j, t, \vec{t}) \in \vec{\mathcal{E}}$, the corresponding overlapping time is denoted by $\sigma_{ij}^{t\vec{t}}$, and for every $(i, j, t, \overleftarrow{t}) \in \overleftarrow{\mathcal{E}}$, the corresponding overlapping time is denoted by $\sigma_{ij}^{t\overleftarrow{t}}$. Our objective is to maximize the sum of overlapping times over all the elements of $\vec{\mathcal{E}}$ and $\overleftarrow{\mathcal{E}}$, i.e.,

$$\text{maximize} \quad \sum_{(i,j,t,\vec{t}) \in \vec{\mathcal{E}}} \sigma_{ij}^{t\vec{t}} + \sum_{(i,j,t,\overleftarrow{t}) \in \overleftarrow{\mathcal{E}}} \sigma_{ij}^{t\overleftarrow{t}}$$

subject to the constraints present in the system. (4)

We model $\sigma_{ij}^{t\vec{t}}$ for all $(i, j, t, \vec{t}) \in \vec{\mathcal{E}}$ and $\sigma_{ij}^{t\overleftarrow{t}}$ for all $(i, j, t, \overleftarrow{t}) \in \overleftarrow{\mathcal{E}}$ in terms of the arrival and departure times of trains. Consider the case, when $(i, j, t, \vec{t}) \in \vec{\mathcal{E}}$. We need to ensure that after we apply the optimization strategy, $\vec{t}$ still stays the temporally closest train to the right of $t$. Otherwise, the only way to achieve a positive overlapping time is to synchronize the braking phase of $t$ with the accelerating phase of $\vec{t}$, which might result in a large deviation of event times compared to the original timetable, especially when there is no or very little overlapping to begin with. We write this constraint as follows:

$$\frac{\left(a_j^{\vec{t}} + d_j^{\vec{t}} - a_i^t - d_i^t\right)}{\left(\bar{a}_j^{\vec{t}} + \bar{d}_j^{\vec{t}} - \bar{a}_i^t - \bar{d}_i^t + \epsilon\right)} \geq 0. \tag{5}$$

Here $\epsilon$ is a very small positive number to prevent division by zero.

Let us denote the start of the braking phase of train $t$ before arriving at platform $i$ by $a_i^{t-}$ and the end of its accelerating phase after departing from the same platform by $d_i^{t+}$. For all $t \in \mathcal{T}$ and for all $i \in \mathcal{N}^t$, the duration of the associated braking phase is $\beta_i^t = a_i^t - a_i^{t-}$ and the duration of the associated accelerating phase is $\alpha_i^t = d_i^{t+} - d_i^t$. Depending on the trip time of the associated trip in consideration, both the durations are within some time bounds, i.e; $\alpha_i^t \in [\underline{\alpha}_i^t, \overline{\alpha}_i^t]$ and $\beta_i^t \in [\underline{\beta}_i^t, \overline{\beta}_i^t]$. Though we do not know the optimal trip time of the trains in advance, these lower and upper bounds can be calculated by existing software [19, page 3]. For the same reason, the start of the braking phase and end of the accelerating phase are within time bounds described by $a_i^{t-} \in [a_i^t - \overline{\beta}_i^t, a_i^t - \underline{\beta}_i^t]$ and $d_i^{t+} \in [d_i^t + \underline{\alpha}_i^t, d_i^t + \overline{\alpha}_i^t]$. All the time bounds are on the order of seconds, as the trip time variation are on the order of seconds, so it is reasonable to pursue a robust formulation. To model the overlapping time $\sigma_{ij}^{t\vec{t}}$ for all $(i, j, t, \vec{t}) \in \vec{\mathcal{E}}$, we propose the following lemma.

**Lemma 1.** *For all $(i, j, t, \vec{t}) \in \vec{\mathcal{E}}$, the overlapping time $\sigma_{ij}^{t\vec{t}}$ between the accelerating phase of $t$ on platform $i$ and the braking phase of $\vec{t}$ on platform $j$, where $(i, j) \in \Omega$, can be modelled by the following robust formulation*

$$a_j^{\vec{t}} - d_i^t + \epsilon \leq \underline{\alpha}_i^t + \underline{\beta}_j^{\vec{t}} + M(1 - \lambda_{ij}^{t\vec{t}}), \tag{6}$$

$$d_i^t - a_j^{\vec{t}} + \epsilon \leq M(1 - \lambda_{ij}^{t\vec{t}}), \tag{7}$$

$$\sigma_{ij}^{t\vec{t}} \geq 0, \tag{8}$$

$$\sigma_{ij}^{t\vec{t}} \leq \underline{\alpha}_i^t \lambda_{ij}^{t\vec{t}}, \tag{9}$$

$$\sigma_{ij}^{t\vec{t}} \leq \underline{\beta}_j^{\vec{t}} \lambda_{ij}^{t\vec{t}}, \tag{10}$$

$$\sigma_{ij}^{t\vec{t}} \leq d_i^t - a_j^{\vec{t}} + \underline{\alpha}_i^t + \underline{\beta}_j^{\vec{t}} + M(1 - \lambda_{ij}^{t\vec{t}}), \tag{11}$$

$$\sigma_{ij}^{t\vec{t}} \leq a_j^{\vec{t}} - d_i^t + M(1 - \lambda_{ij}^{t\vec{t}}). \tag{12}$$

*where $M$ is a large positive number, $\epsilon$ is a small positive number smaller than time granularity considered and $\lambda_{ij}^{t\vec{t}}$ is a binary variable which is one if and only if $\sigma_{ij}^{t\vec{t}}$ is positive.*

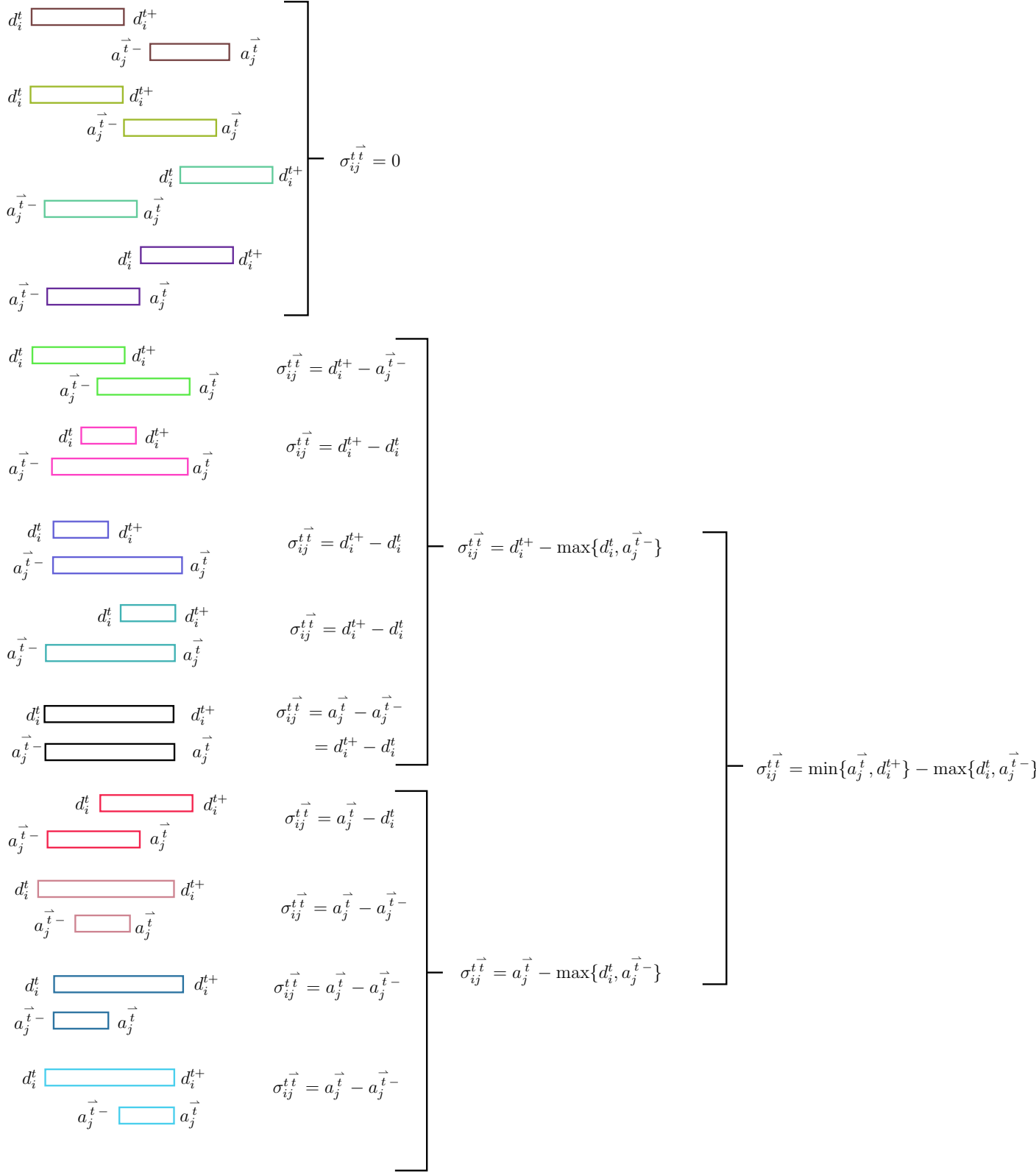

**Figure 2:** All possible overlapping times between the accelerating phase of train $t$ and braking phase of its temporally closest train

*Proof.* We use the hypograph approach to model the overlapping time $\sigma_{ij}^{t\vec{t}}$ in terms of the associated event times [20, page 75, 134]. Consider any $a_i^{t-} \in [a_i^t - \overline{\beta}_i^t, a_i^t - \underline{\beta}_i^t]$ and any $d_i^{t+} \in [d_i^t + \underline{\alpha}_i^t, d_i^t + \overline{\alpha}_i^t]$ for some trip time. From Interval algebra [21], we know that there can be thirteen different kinds of overlapping possible between the accelerating phase of train $t$ and the braking phase of train $\vec{t}$ as shown in Figure 2. However, there is structure among these thirteen relationships when we are concerned with the associated overlapping time $\sigma_{ij}^{t\vec{t}}$. Among them, when $a_j^{\vec{t}-} \geq d_i^{t+}$ or $a_j^{\vec{t}} \leq d_i^t$ there will be no overlapping, i.e., $\sigma_{ij}^{t\vec{t}} = 0$. We will model this scenario using the binary indicator variable $\lambda_{ij}^{tt'}$. If any of the conditions $a_j^{\vec{t}-} \geq d_i^{t+}$ or $a_j^{\vec{t}} \leq d_i^t$ occurs, then $\lambda_{ij}^{t\vec{t}}$ is zero, i.e.,

$$\left((a_j^{\vec{t}-} \geq d_i^{t+}) \vee (a_j^{\vec{t}} \leq d_i^t)\right) \Rightarrow (\lambda_{ij}^{t\vec{t}} = 0) \tag{13}$$

If $\lambda_{ij}^{t\vec{t}}$ is zero, then the overlapping time will be zero, i.e.,

$$(\lambda_{ij}^{t\vec{t}} = 0) \Rightarrow (\sigma_{ij}^{t\vec{t}} = 0)$$
$$\Leftrightarrow (\sigma_{ij}^{t\vec{t}} > 0) \Rightarrow (\lambda_{ij}^{t\vec{t}} = 1), \qquad [\text{as } \sigma_{ij}^{t\vec{t}} \not< 0]. \tag{14}$$

The Equations (13)-(14) are not in mathematical programming format. Using integer programming modelling rules [22, pages 166, 172-174, 183-184], we can model them as follows:

$$a_j^{\vec{t}-} - d_i^{t+} + \epsilon \leq M(1 - \lambda_{ij}^{t\vec{t}}), \tag{15}$$
$$d_i^t - a_j^{\vec{t}} + \epsilon \leq M(1 - \lambda_{ij}^{t\vec{t}}), \tag{16}$$
$$\sigma_{ij}^{t\vec{t}} \geq 0, \tag{17}$$
$$\sigma_{ij}^{t\vec{t}} \leq M\lambda_{ij}^{t\vec{t}}. \tag{18}$$

When $(a_j^{\vec{t}-} < d_i^{t+})$ and $(a_j^{\vec{t}} > d_i^t)$, we can see from Figure 2 that the overlapping time is $\min\{d_i^{t+}, a_j^{\vec{t}}\} - \max\{d_i^t, a_j^{\vec{t}-}\} > 0$. Now we want to model the following logical condition

$$(\lambda_{ij}^{t\vec{t}} = 1) \Rightarrow (\sigma_{ij}^{t\vec{t}} = \min\{d_i^{t+}, a_j^{\vec{t}}\} - \max\{d_i^t, a_j^{\vec{t}-}\})$$

using mathematical programming. Using integer programming modelling rules we can model this situation as follows:

$$\sigma_{ij}^{t\vec{t}} \leq \alpha_i^t \lambda_{ij}^{t\vec{t}}, \tag{19}$$
$$\sigma_{ij}^{t\vec{t}} \leq \beta_j^{\vec{t}} \lambda_{ij}^{t\vec{t}}, \tag{20}$$
$$\sigma_{ij}^{t\vec{t}} \leq d_i^{t+} - a_j^{\vec{t}-} + M(1 - \lambda_{ij}^{t\vec{t}}), \tag{21}$$

$$\sigma_{ij}^{t\vec{t}} \leq a_j^{\vec{t}} - d_i^t + M(1 - \lambda_{ij}^{t\vec{t}}). \tag{22}$$

Combining Equations (15)-(18) and Equations (19)-(22), and using the fact that $M \gg \max\{\alpha_i^t, \beta_j^{\vec{t}}\}$, for any $a_i^{t-} \in [a_i^t - \overline{\beta}_i^t, a_i^t - \underline{\beta}_i^t]$ and any $d_i^{t+} \in [d_i^t + \underline{\alpha}_i^t, d_i^t + \overline{\alpha}_i^t]$ we arrive at the following set of equations:

$$a_j^{\vec{t}-} - d_i^{t+} + \epsilon \leq M(1 - \lambda_{ij}^{t\vec{t}}),$$
$$d_i^t - a_j^{\vec{t}} + \epsilon \leq M(1 - \lambda_{ij}^{t\vec{t}}),$$
$$\sigma_{ij}^{t\vec{t}} \geq 0,$$
$$\sigma_{ij}^{t\vec{t}} \leq \alpha_i^t \lambda_{ij}^{t\vec{t}},$$
$$\sigma_{ij}^{t\vec{t}} \leq \beta_j^{\vec{t}} \lambda_{ij}^{t\vec{t}},$$
$$\sigma_{ij}^{t\vec{t}} \leq d_i^{t+} - a_j^{\vec{t}-} + M(1 - \lambda_{ij}^{t\vec{t}}),$$
$$\sigma_{ij}^{t\vec{t}} \leq a_j^{\vec{t}} - d_i^t + M(1 - \lambda_{ij}^{t\vec{t}}).$$

Using $a_i^{t-} \in [a_i^t - \overline{\beta}_i^t, a_i^t - \underline{\beta}_i^t]$ and $d_i^{t+} \in [d_i^t + \underline{\alpha}_i^t, d_i^t + \overline{\alpha}_i^t]$, the equations above can be transformed into the following robust formulation:

$$a_j^{\vec{t}} - d_i^t + \epsilon \leq \underline{\alpha}_i^t + \underline{\beta}_j^{\vec{t}} + M(1 - \lambda_{ij}^{t\vec{t}}),$$
$$d_i^t - a_j^{\vec{t}} + \epsilon \leq M(1 - \lambda_{ij}^{t\vec{t}}),$$
$$\sigma_{ij}^{t\vec{t}} \geq 0,$$
$$\sigma_{ij}^{t\vec{t}} \leq \underline{\alpha}_i^t \lambda_{ij}^{t\vec{t}},$$
$$\sigma_{ij}^{t\vec{t}} \leq \underline{\beta}_j^{\vec{t}} \lambda_{ij}^{t\vec{t}},$$
$$\sigma_{ij}^{t\vec{t}} \leq d_i^t - a_j^{\vec{t}} + \underline{\alpha}_i^t + \underline{\beta}_j^{\vec{t}} + M(1 - \lambda_{ij}^{t\vec{t}}),$$
$$\sigma_{ij}^{t\vec{t}} \leq a_j^{\vec{t}} - d_i^t + M(1 - \lambda_{ij}^{t\vec{t}}),$$

□

Now consider the case when $(i, j, t, \overleftarrow{t}) \in \overleftarrow{\mathcal{E}}$. Like the previous case, we need to ensure that after we apply the optimization strategy, $\overleftarrow{t}$ still stays the temporally closest train to the left of $t$. An analogous constraint to that of Equation (5) can be easily found by replacing $t$ and $\vec{t}$ in Equation (5) with $\overleftarrow{t}$ and $t$ respectively as follows:

$$\frac{\left(a_i^t + d_i^t - a_j^{\overleftarrow{t}} - d_j^{\overleftarrow{t}}\right)}{\left(\bar{a}_i^t + \bar{d}_i^t - \bar{a}_j^{\overleftarrow{t}} - \bar{d}_j^{\overleftarrow{t}}\right)} \geq 0. \tag{23}$$

Note that the denominator can never be zero on the left hand side of the equation above because of the definition of $\overleftarrow{t}$ in Equation (2). To model the overlapping time $\sigma_{ij}^{t\overleftarrow{t}}$ for all $(i, j, t, \overleftarrow{t}) \in \overleftarrow{\mathcal{E}}$, we propose the following lemma.

**Lemma 2.** *For all $(i, j, t, \overleftarrow{t}) \in \overleftarrow{\mathcal{E}}$, the overlapping time $\sigma_{ij}^{t\overleftarrow{t}}$ between the braking phase of $t$ on platform $i$ and the accelerating phase of $\overleftarrow{t}$ on platform $j$, where $(i, j) \in \Omega$, can be modelled by the following robust formulation*

$$a_i^t - d_j^{\overleftarrow{t}} + \epsilon \leq \underline{\alpha}_i^{\overleftarrow{t}} + \underline{\beta}_j^t + M(1 - \lambda_{ij}^{t\overleftarrow{t}}), \tag{24}$$

$$d_j^{\overleftarrow{t}} - a_i^t + \epsilon \leq M(1 - \lambda_{ij}^{t\overleftarrow{t}}), \tag{25}$$

$$\sigma_{ij}^{t\overleftarrow{t}} \geq 0, \tag{26}$$

$$\sigma_{ij}^{t\overleftarrow{t}} \leq \underline{\alpha}_j^{\overleftarrow{t}} \lambda_{ij}^{t\overleftarrow{t}}, \tag{27}$$

$$\sigma_{ij}^{t\overleftarrow{t}} \leq \underline{\beta}_i^t \lambda_{ij}^{t\overleftarrow{t}}, \tag{28}$$

$$\sigma_{ij}^{t\overleftarrow{t}} \leq d_j^{\overleftarrow{t}} - a_i^t + \underline{\alpha}_i^{\overleftarrow{t}} + \underline{\beta}_j^t + M(1 - \lambda_{ij}^{t\overleftarrow{t}}), \tag{29}$$

$$\sigma_{ij}^{t\overleftarrow{t}} \leq a_i^t - d_j^{\overleftarrow{t}} + M(1 - \lambda_{ij}^{t\overleftarrow{t}}), \tag{30}$$

*where $M$ is a large positive number, $\epsilon$ is a small positive number smaller than time granularity considered and $\lambda_{ij}^{t\overleftarrow{t}}$ is a binary variable which is one if and only if $\sigma_{ij}^{t\overleftarrow{t}}$ is positive.*

*Proof.* The lemma can be easily proved by replacing $i$, $j$, $t$ and $\overrightarrow{t}$ with $j$, $i$, $\overleftarrow{t}$ and $t$ respectively in Lemma 1. □

## 4 Modelling the constraints

The constraints in the railway network show how the events are related. In this section we describe, model and justify the constraints.

### 4.1 Trip time constraint

Consider the trip of any train $t \in \mathcal{T}$ from platform $i$ to platform $j$ along the track $(i, j) \in \mathcal{A}^t$. The train $t$ departs from platform $i$ at time $d_i^t$ and arrives at platform $j$ at time $a_j^t$ and the train can have a trip time between $\underline{\tau}_{ij}^t$ and $\overline{\tau}_{ij}^t$. The trip time constraint can be written as follows:

$$\forall t \in \mathcal{T} \; \forall (i, j) \in \mathcal{A}^t \quad \left( \underline{\tau}_{ij}^t \leq a_j^t - d_i^t \leq \overline{\tau}_{ij}^t \right). \tag{31}$$

## 4.2 Dwell time constraint

When any train $t \in \mathcal{T}$ arrives at a platform $i \in \mathcal{N}^t$, it dwells there for a certain time interval denoted by $[\underline{\delta}_i^t, \overline{\delta}_i^t]$ so that the passengers can get off and get on the train. After the dwelling time is over, the train departs from the station. The difference between the departure time $d_i^t$ and arrival time $a_i^t$ corresponding to the dwelling mentioned lies between $\underline{\delta}_i^t$ and $\overline{\delta}_i^t$. The dwell time constraint can be written as follows:

$$\forall t \in \mathcal{T}\ \forall i \in \mathcal{N}^t \quad (\underline{\delta}_i^t \leq d_i^t - a_i^t \leq \overline{\delta}_i^t). \tag{32}$$

Every train $t \in \mathcal{T}$ arrives at the first platform $\mathcal{N}^t(1)$ in its train-path either from the depot or by turning around from some other line, and departs from the final platform $\mathcal{N}^t(|\mathcal{N}^t|)$ in order to either return to the depot or start as a new train on another line by turning around. So, the train $t$ dwells at all the platforms in $\mathcal{N}^t$. This is the reason why in Equation (32) the platform index $i$ is varied over all the elements of the set $\mathcal{N}^t$.

## 4.3 Connection constraint:

In many cases, a single train might not exist between the origin and desired destination of a passenger. To circumvent such issues, connecting trains are often used by the railway management at interchange stations. Let $\chi \subseteq \mathcal{N} \times \mathcal{N}$ be the set of platform pairs where passengers transfer between trains. If $(i, j) \in \chi$, then both the platforms $i$ and $j$ are situated at the same station, and there exists a train $t \in \mathcal{T}$ arriving at platform $i$ and another train $t' \in \mathcal{T}$ departing from platform $j$ such that a connection time window needs to be maintained between train $t$ and $t'$ for passengers to get off from the first train and get on the latter. Note that order matters here. Let $\mathcal{C}_{ij}$ be the set of connecting or turning-around train pairs for a platform pair $(i, j) \in \chi$. Then the connection constraint can be written as:

$$\forall (i, j) \in \chi\ \forall (t, t') \in \mathcal{C}_{ij} \quad (\underline{\chi}_{ij}^{tt'} \leq d_j^{t'} - a_i^t \leq \overline{\chi}_{ij}^{tt'}), \tag{33}$$

where $\underline{\chi}_{ij}^{tt'}$ and $\overline{\chi}_{ij}^{tt'}$ are the lower bound and upper bound of the time window to achieve the described connection between the associated trains.

## 4.4 Turn-around constraint:

Turn-around events provide a connection between planned terminal-to-terminal trips that will be executed by a single physical train. This constraint is analogous to the trip time constraint for these connecting movements. After a train arrives at the terminal platform of a train-line, it often turns around by traversing the crossing-over and starts travelling through another train-line. From a railway management perspective, a time window has to be maintained between the departure of the train from the terminal platform of the first line and the arrival time of it (labelled as a different train, though it is physically the same train) on another platform of the second line. Let $\varphi$ be the set of all crossing-overs, where turn-around events occur. Consider any crossing-over $(i, j) \in \varphi$, where the platforms $i$ and $j$ are situated on different train-lines. Let $\mathcal{B}_{ij}$ be the set of all train pairs involved in corresponding turn-around events. Let $(t,t') \in \mathcal{B}_{ij}$. Train $t \in \mathcal{T}$ turns around at platform $i$ by travelling through the crossing-over $(i,j)$, and from platform $j$ starts traversing a different train-line as train $t' \in \mathcal{T} \setminus \{t\}$. A time window denoted by $[\underline{\kappa}_{ij}^{tt'}, \overline{\kappa}_{ij}^{tt'}]$ has to be maintained between the mentioned events, where $\underline{\kappa}_{ij}^{tt'}$ and $\overline{\kappa}_{ij}^{tt'}$ are the lower bound and upper bound respectively. We can write this constraint as follows:

$$\forall (i,j) \in \varphi \; \forall (t,t') \in \mathcal{B}_{ij} \quad (\underline{\kappa}_{ij}^{tt'} \leq a_j^{t'} - d_i^t \leq \overline{\kappa}_{ij}^{tt'}). \tag{34}$$

## 4.5 Headway constraint:

In any railway network, a minimum amount of time between the departures of consecutive trains is always maintained. This time is called headway time. Let $(i,j) \in \mathcal{A}$ be the track between two platform $i$ and $j$, and $\mathcal{H}_{ij}$ be the set of train-pairs who move along that track successively on the order of their departures. Now, assume train $t$ and train $t'$ move along this track in same direction where $(t,t') \in \mathcal{H}_{ij}$. Let $h_i^{tt'}$ and $h_j^{tt'}$ be the associated headway times at platform $i$ and platform $j$ respectively. So, the headway constraint can be written as:

$$\forall (i,j) \in \mathcal{A} \; \forall (t,t') \in \mathcal{H}_{ij}$$
$$(h_i^{tt'} \leq d_i^{t'} - d_i^t \; \wedge \; h_j^{tt'} \leq d_j^{t'} - d_j^t). \tag{35}$$

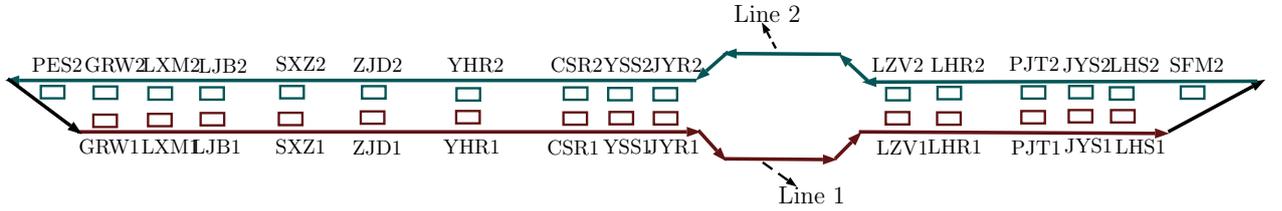

**Figure 3:** Railway network considered for numerical study

## 4.6 Total travel time constraint:

To maintain the quality of service in the railway network, it is desired that for every train $t \in \mathcal{T}$, the total travel time to traverse its train-path stays within a time window $[\underline{\tau}_\mathcal{P}^t, \overline{\tau}_\mathcal{P}^t]$, where $\underline{\tau}_\mathcal{P}^t$ and $\overline{\tau}_\mathcal{P}^t$ are the corresponding lower and upper bound respectively. We can write this constraint as follows:

$$\forall t \in \mathcal{T} \quad (\underline{\tau}_\mathcal{P}^t \leq a_{\mathcal{N}^t(|\mathcal{N}^t|)}^t - d_{\mathcal{N}^t(1)}^t \leq \overline{\tau}_\mathcal{P}^t). \tag{36}$$

Here $\mathcal{N}^t(1)$ and $\mathcal{N}^t(|\mathcal{N}^t|)$ are the first and last platform in the train-path of $t$.

# 5 Full optimization model

In this section we collect the objective and all the constraints discussed in the previous two sections, and propose our optimization problem to maximize the total duration of overlapping times of the SPSTPs in order to utilize regenerative braking energy produced by trains in a railway network. The full optimization model is as follows:

$$\text{maximize} \sum_{(i,j,t,\overrightarrow{t})\in \overrightarrow{\mathcal{E}}} \sigma_{ij}^{t\overrightarrow{t}} + \sum_{(i,j,t,\overleftarrow{t})\in \overleftarrow{\mathcal{E}}} \sigma_{ij}^{t\overleftarrow{t}}$$

subject to

Equations (31), (32), (33), (34), (35) and (36)

$\forall (i,j,t,\overrightarrow{t}) \in \overrightarrow{\mathcal{E}}$    Equations (5),(6)-(12)

$\forall (i,j,t,\overleftarrow{t}) \in \overleftarrow{\mathcal{E}}$    Equations (23),(24)-(30)

$\forall t \in \mathcal{T} \; \forall i \in \mathcal{N}^t \quad (a_i^t \geq 0, d_i^t \geq 0)$

$\forall (i,j,t,\tilde{t}) \in \mathcal{E} \quad (\lambda_{ij}^{t\tilde{t}} \in \{0,1\}, \sigma_{ij}^{t\tilde{t}} \geq 0),$

where the decision variables are $a_i^t$, $d_i^t$, $\lambda_{ij}^{t\tilde{t}}$ and $\sigma_{ij}^{t\tilde{t}}$.

As the model is a MIP with bounds, the optimization problem is $\mathcal{NP}$-hard [23, page 242]. However, in the next section we show that for the size of the railway data considered in practice, the running time is quite acceptable.

# 6 Numerical experiment

In this section we apply our model to nine different problem instances to service PES2-SFM2 of line 8 of Shanghai Metro network. The numerical study was executed on a Intel Core i7-46400 CPU with 8 GB RAM running Windows 8.1 Pro operating system. To model our problem, we have used `JuMP` - an open source algebraic modelling language embedded in the programming language `Julia` [24]. Within our `JuMP` code we have called academic version of Gurobi Optimizer 6.0 as the solver. In each instance we have an initial feasible timetable with a duration of six hours. In most of the railway networks the duration of the off-peak or rush hours is smaller than six hours, so a timetable that six hours is sufficient for practical purpose. The number of trains increases as the average headway time decreases. The results of the numerical study are shown in Table 1.

The service PES2-SFM2 of line 8 of the Shanghai Metro network is shown in Figure 3. There are two lines in this network: Line 1 and Line 2. There are fourteen stations in the network denoted by all capitalized words in the figure. Each station has two platforms each on different train lines, *e.g.*, LXM is station that has two opposite platforms: LXM1 and LXM2 on Line 1 and Line 2 respectively. The platforms are denoted by rectangles. The platforms indicated by PES2 and SFM2 are the turn-around points on Line 2, and the crossing-overs are PES2-GRW1 and LHS1-SFM2.

**Table 1:** Results of the numerical study (running time 1200 s)

| Number of Trains | Number of Constraints | Binary Variables | Continuous Variables | CPU Time (s) | Nodes Explored | Intial Overlapping Time (s) | Final Overlapping Time (s) | Initial Effective Energy Consumption (kWh) | Final Effective Energy Consumption (kWh) | Reduction in Effective Energy Consumption (%) |
|---|---|---|---|---|---|---|---|---|---|---|
| 100 | 10437 | 183 | 3133 | 86.64 | 4534 | 454.26 | 2828.09 | 8594.27 | 7698.28 | 10.42 |
| 112 | 12166 | 272 | 3576 | 7.26 | 459 | 1307.23 | 4354.5 | 14500.05 | 13740.82 | 5.23 |
| 124 | 12684 | 188 | 3846 | 0.28 | 0 | 1190.89 | 3119.72 | 8524.96 | 7618.28 | 10.63 |
| 134 | 13186 | 128 | 4081 | 0.25 | 0 | 315.85 | 1950.88 | 5451.29 | 4637.33 | 14.93 |
| 146 | 16378 | 426 | 4733 | 0.45 | 0 | 2412.92 | 7314.48 | 20223.63 | 16998.21 | 15.95 |
| 156 | 16264 | 278 | 4880 | 15.34 | 6 | 2096.23 | 4226.93 | 14324.78 | 13319.2 | 7.02 |
| 168 | 18154 | 390 | 5346 | 56.94 | 1164 | 1823.1 | 6299.66 | 20274.44 | 17817.97 | 12.12 |
| 178 | 19650 | 472 | 5723 | 0.51 | 0 | 2874.67 | 8195.41 | 21685.04 | 17619.35 | 18.75 |
| 190 | 22947 | 785 | 6390 | 9.6 | 0 | 5858.58 | 12924.23 | 35613.96 | 33574.48 | 5.73 |

The feasible timetables are provided to us by Thales Canada Inc.. We have taken $M = 1000$ and $\epsilon = 0.005$. We have applied our optimization model to find the optimal timetable that maximizes the total overlapping time of the SPSTPs. We see from Table 1 that for each instance, our optimization model produces an optimal timetable with significant increase in the total overlapping time in comparison with the initial timetable. Such increase in the total overlapping time would make it possible to save significant amount of electrical energy produced by the braking trains by transferring it to the accelerating trains via the overhead contact lines. We can see that, in all of the cases our model has found the optimal timetables very quickly, the largest runtime being 86.64s. Number of nodes explored is quite small for every instance, and for five out of nine instances an optimal solution is found at the root node. After we get the final timetable, we calculate the total *effective energy consumption* by all the trains involved in SPSTPs and compare it with the original timetables. The effective energy consumption of a train for a trip is defined as the difference between total energy required to make a trip and the amount of energy that is being supplied by a braking train during synchronization process. So, the effective energy consumption is the energy that will be consumed from the electrical substations.

The energy calculation is done using `SPSIM` [19], and `Cubature` [25]. `SPSIM` is a proprietary software owned by Thales Canada Inc. that calculates the power versus time graphs of all the active trains for the original and optimal timetables. Cubature is an open-source `Julia` package written by Steven G. Johnson, which is used to calculate the effective area under the power versus time graphs to determine 1) the total energy required by the trains during the trips, 2) the total transferred regenerative energy during the SPSTPs, and 3) the effective energy consumption that the difference of the first two quantities. The effective energy consumption of the optimal timetables in comparison with the original ones gets reduced for every instance, with smallest reduction being 5.23% and the largest reduction

being 18.74%.

# 7 Conclusion

In this paper we have presented a robust mixed integer optimization model to utilize regenerative braking energy produced by trains in a railway network. The optimization model calculates a railway timetable that saves regenerative energy of braking trains by transferring it to suitable accelerating trains in need of energy. We have used the hypograph approach and interval algebra to formulate the objective function. We have modelled different constraints to describe a feasible timetable for a railway network of arbitrary topology. We have applied our optimization model to different instances of a railway network for a time horizon spanning six hours. For each instance, our model has found an optimal timetable very quickly with significant reduction in the effective energy consumption by the trains.